# CRITICAL THRESHOLDS IN 2D RESTRICTED EULER-POISSON EQUATIONS

HAILIANG LIU AND EITAN TADMOR

Contents

ABSTRACT. We provide a complete description of the critical threshold phenomena for the two-dimensional localized Euler-Poisson equations, introduced by the authors in [17]. Here, the questions of global regularity vs. finite-time breakdown for the 2D Restricted Euler-Poisson solutions are classified in terms of precise explicit formulae, describing a remarkable variety of critical threshold surfaces of initial configurations. In particular, it is shown that the 2D critical thresholds depend on the relative size of three quantities: the initial density, the initial divergence as well as the initial spectral gap, that is, the difference between the two eigenvalues of the $2 \times 2$ initial velocity gradient.

## 1. INTRODUCTION AND STATEMENT OF MAIN RESULTS

We are concerned with the critical threshold phenomena in multi-dimensional Euler-Poisson equations. In this paper we consider a localized version of the following two-dimensional Euler-Poisson equations

(1.1) $\quad\quad\quad\quad \partial_t \rho + \nabla \cdot (\rho U) = 0, \quad x \in \mathbb{R}^2, \quad t \in \mathbb{R}^+,$

(1.2) $\quad\quad\quad\quad \partial_t(\rho U) + \nabla \cdot (\rho U \otimes U) = -k\rho \nabla \phi,$

(1.3) $\quad\quad\quad\quad -\Delta \phi = \rho - c, \quad x \in \mathbb{R}^2,$

which are the usual statements of the conservation of mass, Newton's second law, and the Poisson equation defining, say, the electric field in terms of the charge. Here $k > 0$ is a scaled physical constant, which signifies the property of the underlying repulsive forcing (avoiding the case of an attractive force with $k < 0$), and $c$ denotes the constant "background" state. The unknowns are the local density $\rho = \rho(x,t)$, the velocity field $U = (u, v)(x, t)$, and the potential $\phi = \phi(x, t)$. It follows that as long as the solution remains smooth, the velocity $U$ solves a forced transport equation

$$\partial_t U + U \cdot \nabla U = F, \quad F = -k\nabla\phi \tag{1.4}$$

with $\phi$ being governed by Poisson's equation (1.3).

This hyperbolic-elliptic coupled system (1.1)-(1.3) describes the dynamic behavior of many important physical flows including the charge transport [22], plasma with collision [13], cosmological waves [2] and the expansion of the cold ions [12]. Let us mention that the Euler-Poisson equations could also be realized as the semi-classical limit of Schrödinger-Poisson equation and are found in the 'cross-section' of Vlasov-Poisson equations. These relations have been the subject of a considerable amount of work in recent years, and we refer to [10], [6] and references therein for further details.

To put our study in a proper perspective we recall that there has been a considerable amount of literature available on the global behavior of Euler-Poisson and related problems, from local existence in the small $H^s$-neighborhood of a steady state [18, 23, 9] to global existence of weak solutions with geometrical symmetry [5], for the two-carrier types in one dimension [29], the relaxation limit for the weak entropy solution, consult [21] for isentropic case, and [14] for isothermal case.

For the question of global behavior of strong solutions, however, the choice of the initial data and/or damping forces is decisive. The non-existence results in the case of attractive forces, $k < 0$, have been obtained by Makino-Perthame [20], and for repulsive forces by Perthame [24]. For the study on the singularity formation in the model with diffusion and relaxation, consult [30]. In all these cases, the finite life span is due to a *global* condition of large enough initial (generalized) energy, staying outside a critical threshold ball. Using the characteristic-based method, Engelberg [7] gave local conditions for the finite-time loss of smoothness of solutions in Euler-Poisson equations. Global existence due to damping relaxation and with non-zero background can be found in [27, 28, 15]. For the model without damping relaxation the global existence was obtained by Guo [11] assuming the flow is irrotational. His result applies to $H_2$-small neighborhood of constant state.

When dealing with the questions of time regularity for Euler-Poisson equations without damping, one encounters several limitations with the classical stability analysis. Among others issues, we mention that

(i) the stability analysis does not tell us how large perturbations are allowed before losing stability – indeed, the smallness of the initial perturbation is essential to make the energy method work, e.g., [11];

(ii) the steady solution may be only conditionally stable due to the weak dissipation in the system, say in the 1D Euler-Poisson equations [8].

In order to address these difficulties we advocated, in [8], a new notion of critical threshold (CT), which describes the conditional stability of the 1D Euler-Poisson equations, where the answer to the question of global vs local existence depends on whether the initial configuration crosses an intrinsic, $O(1)$ critical threshold. Little or no attention has been paid to this remarkable phenomena, and our goal is to bridge the gap of previous studies on the behavior in Euler-Poisson solutions, a gap between the regularity of Euler-Poisson solutions in the small and their finite-time breakdown in the large. The critical threshold (CT) in the 1D Euler-Poisson system was completely characterized in terms of the relative size of the initial velocity slope and the initial density. Moving to the multi-D setup, one has first to identify what are the proper quantities which govern the critical threshold phenomena . In [17] we have shown that these quantities depend in an essential manner on the *eigenvalues* of the gradient velocity matrix, $\nabla u$. In order to trace the evolution of $M := \nabla U$, we differentiate (1.4), obtaining formally

(1.5) $$\partial_t M + U \cdot \nabla M + M^2 = -k(\nabla \otimes \nabla)\phi = kR[\rho - c],$$

where $R[\ ]$ is the $2 \times 2$ Risez matrix operator, defined as

$$R[f] =: \nabla \otimes \nabla \Delta^{-1}[f] = \mathcal{F}^{-1}\left\{\frac{\xi_j \xi_k}{|\xi|^2}\hat{f}(\xi)\right\}_{j,k=1,2}.$$

The above system is complemented by its coupling with the density $\rho$ which is governed

(1.6) $$\partial_t \rho + U \cdot \nabla \rho + \rho tr M = 0.$$

Passing to the Lagrangian coordinates, that is, using the change of variables $\alpha \mapsto x(\alpha, t)$ with $x(\alpha, t)$ solving

$$\frac{dx}{dt} = U(x, t), \quad x(\alpha, 0) = \alpha,$$

Euler-Poisson equations are recast into the coupled system

(1.7) $$\frac{d}{dt}M + M^2 = kR[\rho - c],$$

(1.8) $$\frac{d}{dt}\rho + \rho tr M = 0,$$

with $d/dt$ standing for the usual material derivative, $\partial_t + U \cdot \nabla$. It is the global forcing, $kR[\rho - c]$, which presents the main obstacle to study the CT phenomena of the multi-D Euler-Poisson setting.

In this work we focus on the Restricted Euler-Poisson (REP) system introduced in [17], which is obtained from (1.7) by restricting attention to the local isotropic trace, $\frac{k}{2}(\rho - c)I_{2\times 2}$, of the global coupling term $kR[\rho - c]$, namely

(1.9) $$\frac{d}{dt}M + M^2 = \frac{k}{2}(\rho - c) \cdot I_{2\times 2},$$

(1.10) $$\frac{d}{dt}\rho + \rho tr M = 0.$$

We are concerned with the initial value REP problem (1.9),(1.10), subject to initial data

$$(M, \rho)(\cdot, 0) = (M_0, \rho_0).$$

We note in passing that the REP system is to the full Euler-Poisson equations what the Restricted Euler model is to the full Euler equations, consult [26, 3, 1, 4, 17]. The existence of a critical threshold phenomena associated with this 2D REP model with zero background, $c = 0$, was first identified by us [17]. The current paper provides a precise description of the critical threshold for the 2D REP system (1.9), (1.10), with both zero or non-zero background charge. In particular, we use the so called Spectral Dynamics Lemma, [17, Lemma 3.1] to obtain remarkable explicit formulae for the critical threshold surfaces summarized in the main Theorems 1.1 and 1.2 below.

To state our main results, we introduce two quantities with which we characterize the behavior of the velocity gradient tensor $M$. These are the trace, $d := trM$, and we note that in case $M$ coincides with $\nabla U$ then $d$ stand for the divergence, $d = u_x + v_y$; and the nonlinear quantity $\Gamma := (trM)^2 - 4detM$, which serves as an index for the *spectral gap*. Indeed, if $\lambda_i, i = 1, 2$ are the eigenvalues of $M$, then

$$\lambda_1 = \frac{1}{2}[d - \sqrt{\Gamma}], \quad \lambda_2 = \frac{1}{2}[d + \sqrt{\Gamma}]$$

and hence $\Gamma$ is nothing but the square of the spectral gap $\Gamma = (\lambda_2 - \lambda_1)^2$. We note that when $M$ coincides with $\nabla U$, then $\Gamma = (u_x - v_y)^2 + 4u_y v_x$, and the role of this spectral gap was first identified in the context of the 2D Eikonal equation in [17, lemma 5.2].

We observe that if $\Gamma < 0$ then the spectral gap is purely imaginary, and otherwise the spectral gap is real.

**Theorem 1.1.** *[2D REP with zero background]. Consider the 2D repulsive REP system (1.9)-(1.10) with $k > 0$ and with zero background $c = 0$. Then the solution of 2D REP remains smooth for all time if and only if the initial data $(\rho_0, M_0)$ lies in one of the following two regions, $(\rho_0, d_0, \Gamma_0) \in S_1 \cup S_2$:*

*(i) Either*

$$(\rho_0, d_0, \Gamma_0) \in S_1, \quad S_1 := \left\{(\rho, d, \Gamma) \,\Big|\, \Gamma \leq 0 \quad and \quad \begin{cases} d \geq 0 & \text{if } \rho = 0 \\ d \text{ arbitrary} & \text{if } \rho > 0 \end{cases}\right\};$$

*(ii) or*

$$(\rho_0, d_0, \Gamma_0) \in S_2, \quad S_2 := \left\{(\rho, d, \Gamma) \,\Big|\, \rho > 0, \ \Gamma > 0, \quad and \quad d \geq g(\rho, \Gamma)\right\}$$

*where*

$$g(\rho, \Gamma) := \text{sgn}(\Gamma - 2k\rho)\sqrt{\Gamma - 2k\rho + 2k\rho ln\left(\frac{2k\rho}{\Gamma}\right)}.$$

**Theorem 1.2.** *[2D REP with nonzero background]. Consider the 2D repulsive REP system (1.9)-(1.10) with $k > 0$ and with non-zero background $c > 0$. Then the solution of 2D REP remains smooth for all time if and only if the initial data $(\rho_0, M_0)$ lies in one of the following three regions, $(\rho_0, d_0, \Gamma_0) \in S_1 \cup S_2 \cup S_3$:*

*(i) Either*

$$(\rho_0, d_0, \Gamma_0) \in S_1, \quad S_1 := \left\{(\rho, d, \Gamma) \,\Big|\, \Gamma \leq 0 \quad and \quad \begin{cases} d \geq 0 & \text{if } \rho = 0 \\ d \text{ arbitrary} & \text{if } \rho > 0 \end{cases}\right\};$$

*(ii) or*

$$(\rho_0, d_0, \Gamma_0) \in S_2, \quad S_2 := \left\{(\rho, d, \Gamma) \,\Big|\, 0 < \Gamma < \frac{k}{2c}\rho^2 \quad and \quad \begin{cases} |d| \leq g_1(\rho, \Gamma) & \text{if } \Gamma < 2k(\rho - c) \\ d \geq g_1(\rho, \Gamma) & \text{if } \Gamma \geq 2k(\rho - c) \end{cases}\right\},$$

*where*

$$g_1(\rho, \Gamma) := \sqrt{\Gamma - 2k\left[c + \sqrt{\rho^2 - 2ck^{-1}\Gamma} + \rho ln\left(\frac{\rho - \sqrt{\rho^2 - 2ck^{-1}\Gamma}}{2c}\right)\right]};$$

*(iii) or*

$$(\rho_0, d_0, \Gamma_0) \in S_3, \qquad S_3 := \left\{ (\rho, d, \Gamma) \;\Big|\; \Gamma = \frac{k}{2c}\rho^2, \quad d = g_2(\rho, \Gamma), \quad \rho > 0 \right\},$$

*where*

$$g_2(\rho) = g_1(\rho, \Gamma)|_{\Gamma = \frac{k}{2c}\rho^2} := \sqrt{-2ck + \frac{k}{2c}\rho^2 + 2k\rho \ln\left(\frac{2c}{\rho}\right)}.$$

Several remarks are in order.

1. The above results show that the global smooth solution is ensured if the *initial* velocity gradient has complex eigenvalues, which applies, for example, for a class of initial configurations with sufficiently large vorticity $|u_{0y} - v_{0x}| \gg 1$. With other initial configurations, however, the finite time breakdown of solutions may – and actually does occur, unless the initial divergence is above a critical threshold, expressed in terms of the initial density and initial spectral gap. Hence, global regularity depends on whether the initial configuration crosses an intrinsic, $\mathcal{O}(1)$ critical threshold.

2. The critical threshold in the 1D Euler-Poisson equations depends on the relative size of the initial velocity slope and the initial density, consult [8]. In contrast to the 1D scenario, the critical threshold presented here depends on three initial quantities: density $\rho_0$, divergence $\nabla \cdot U_0$ and initial spectral gap $\Gamma_0 = (u_{0x} - v_{0y})^2 + 4u_{0y}v_{0x}$.

3. Theorem 1.1 tells us that the size of initial sub-critical range which gives rise to regular solution is decreasing as the initial ratio $\Gamma_0/\rho_0$ is increasing. In particular when this ratio is larger than $2k$, then the initial divergence must stay above a positive critical threshold to avoid the finite time breakdown.

4. From Theorem 1.2 we see that the initial critical range which guarantees global regularity shrinks as the initial ratio $\Gamma_0/\rho_0^2$ is increasing in $(-\infty, \frac{k}{2c})$. Finite time breakdown must occur when this ratio is larger than $\frac{k}{2c}$.

5. The limit $c \downarrow 0$ is a sort of a *singular limit* and hence one cannot recover Theorem 1.1 simply by passing to the limit $c \to 0$ in Theorem 1.2. □

It is well known that a finite time breakdown is a generic phenomena for nonlinear hyperbolic convection equations, which are realized by the formation of shock discontinuities. In the context of Euler-Poisson equations, however, there is a delicate balance between the forcing mechanism (governed by Poisson equation), and the nonlinear focusing (governed by Newton's second law), which supports a critical threshold phenomena.

In this paper we show how the persistence of the global features of the solutions for REP hinges on a delicate balance between the nonlinear convection and the localized forcing mechanism dictated by the Poisson equation. Here we use these restricted models to demonstrate the ubiquity of *critical thresholds* in the solutions of some of the equations of mathematical physics. This remarkable CT phenomena has been found in other contexts, say convolution model for nonlinear conservation laws [16], a nonlocal model in the nonlinear wave propagation [25], etc.

We note in passing that in this paper we restrict ourselves to Euler-Poisson equations with localized forcing. The presence of the global forcing allows for additional balance, and we hope to explore the critical threshold phenomena for the general model with global forcing in a future work.

We now conclude this section by outlining the rest of the paper. In Section 2 we study the critical threshold for the REP with zero background. The key observation is that the spectral gap is conserved along particle path. With this property we will be able to reduce the full dynamics on the 2D manifold parameterized by this initial spectral gap. In Section 3 we discuss the critical threshold for the REP with non-zero background, where the CT arguments become considerably more involved. We treat the different cases which are indexed by the initial spectral gap.

## 2. 2D Restricted EP with zero background

In this section we prove the existence of the critical threshold of the 2D restricted EP with zero background ($c = 0$)

$$(2.1) \qquad \frac{d}{dt}M + M^2 = \frac{k}{2}\rho I_{2\times 2},$$

$$(2.2) \qquad \frac{d}{dt}\rho + \rho tr M = 0.$$

This system with initial data $(\rho_0, M_0)$ is well-posed in the usual $H^s$ Sobolev spaces for a short time. The global regularity follows from the standard boot-strap argument once a priori estimate on $\|M(\cdot)\|_{L^\infty}$ is obtained. First we show that for the 2D restricted EP (2.1)-(2.2), the velocity gradient tensor is completely controlled by the divergence $d$ and the density $\rho$.

**Lemma 2.1.** *Let $M$ be the solution of the 2D restricted EP, then the bounded of $M$ depends on the boundedness of $trM$ and $\rho$, namely, there exist constant, $Const = Const_T$ such that*

$$\|M(\cdot, t)\|_{L^\infty[0,T]} \leq Const_T . \|(trM, \rho)\|_{L^\infty[0,T]}.$$

*Proof.* For the 2-D case the velocity gradient tensor is completely governed by $p := M_{11} - M_{22}$, $q := M_{12} + M_{21}$, $\omega = M_{12} - M_{21}$ and $d = M_{11} + M_{22}$. From the $M$-equation (2.1),

$$\frac{d}{dt}\begin{pmatrix} M_{11} & M_{12} \\ M_{21} & M_{22} \end{pmatrix} + \begin{pmatrix} M_{11}^2 + M_{21}M_{12} & dM_{12} \\ dM_{21} & M_{21}M_{12} + M_{22}^2 \end{pmatrix} = \frac{k}{2}\rho I_{2\times 2},$$

one can obtain

$$\frac{d}{dt}p + pd = 0,$$
$$\frac{d}{dt}q + qd = 0,$$
$$\frac{d}{dt}\omega + \omega d = 0,$$

which when combined with the mass equation

$$\frac{d}{dt}\rho + \rho d = 0$$

gives

$$(p, q, \omega) = (p_0, q_0, \omega_0)\rho_0^{-1}\rho.$$

This shows that $|M_{ij}|_{L^\infty}$ are bounded in terms of $|d|_{L^\infty}$ and $|\rho|_{L^\infty}$ as asserted. □

This lemma tells us that to show the global regularity it suffices to control the divergence $d$ and the density $\rho$. Let $\lambda_i, i = 1, 2$, be the eigenvalues of the velocity gradient tensor, then $d = \lambda_1 + \lambda_2$ and the continuity equation (1.10) reads

$$\frac{d}{dt}\rho + \rho(\lambda_1 + \lambda_2) = 0. \tag{2.3}$$

The Spectral Dynamics Lemma [17, Lemma 3.1] tells us that the velocity gradient equation (1.9) yields

$$\frac{d}{dt}\lambda_1 + \lambda_1^2 = \frac{k}{2}\rho, \tag{2.4}$$

$$\frac{d}{dt}\lambda_2 + \lambda_2^2 = \frac{k}{2}\rho. \tag{2.5}$$

Following [17] we consider the difference of the last two equations which gives for $\eta := \lambda_2 - \lambda_1$,

$$\frac{d}{dt}\eta + \eta(\lambda_1 + \lambda_2) = 0.$$

This combined with the mass equation (2.3) and $trM = \lambda_1 + \lambda_2$ yields

$$\frac{d}{dt}\left(\frac{\eta}{\rho}\right) = 0 \Rightarrow \frac{\eta}{\rho} = \frac{\eta_0(\alpha)}{\rho_0(\alpha)}, \quad \alpha \in \mathsf{IR}^2.$$

Set $\beta := \eta_0^2(\alpha)/\rho_0^2(\alpha)$ as a moving parameter with the initial position $\alpha \in \mathsf{IR}^2$, one then obtains a closed system for $\rho$ and $d$

$$\rho' + \rho d = 0, \quad ' := \frac{d}{dt}, \tag{2.6}$$

$$d' + \frac{d^2 + \beta\rho^2}{2} = k\rho. \tag{2.7}$$

The first is the mass equation; the second is a re-statement of summing (2.4),(2.5), $d' + (d^2 + \eta^2)/2 = k\rho$ with $\eta^2 = \beta\rho^2$.

We shall study the dynamics of $(\rho, d)$ parameterized by $\beta$. It is easy to see that if the initial eigenvalues are complex, then the eigenvalues remain complex as time evolves. From

$$\beta = \frac{\Gamma_0}{\rho_0^2}, \quad \Gamma_0 = (\lambda_2(0) - \lambda_1(0))^2$$

we see that we need to distinguish between two cases, namely, $\beta < 0$ where the initial spectral gap is complex, and $\beta \geq 0$ where the initial spectral gap is real.

2.1. **Complex Spectral Gap.** We first study the case $\beta < 0$ when the initial eigenvalues are complex, i.e. $Im(\lambda_i) \neq 0$.

**Lemma 2.2.** *The solution of 2-D REP remains smooth for all time if eigenvalues are initially complex. Moreover there is a global invariant given by*

$$\frac{d^2 - \beta\rho^2}{\rho} + 2kln\rho = Const. \tag{2.8}$$

*Proof.* To obtain the desired global invariant we set $q := d^2$, then from the equations (2.6)-(2.7) we deduce

$$\frac{dq}{d\rho} = 2d\frac{d'}{\rho'} = -2k + \beta\rho + \frac{q}{\rho}.$$

Integration gives

$$\frac{q}{\rho} - \beta\rho + 2kln\rho = Const,$$

which leads to (2.8). The boundedness of $d$ follows at once, since for negative $\beta$'s,

$$d^2 \leq \max_{\rho>0}\left\{Const.\rho - 2kln\rho + \beta\rho^2\right\} =: C_1^2.$$

In particular, substitution of the lower bound, $d \geq -C_1$ into the mass equation gives

$$\rho' \leq C_1\rho$$

which yields the desired upper-bound for the density, $\rho(\cdot,t) \leq \rho_0(\alpha)e^{C_1 t}$. $\square$

*Remark.* More precise information about the large time behavior is available from phase plane analysis. According to (2.7), the zero level set of $d' = 0$ is the hyperbola $Q := k\rho - (d^2 + \beta\rho^2)/2 = 0$, with a right branch passing critical point $(0,0)$ and a left branch is located in the left half plane, $\rho < 0$, see Figure 2.1.

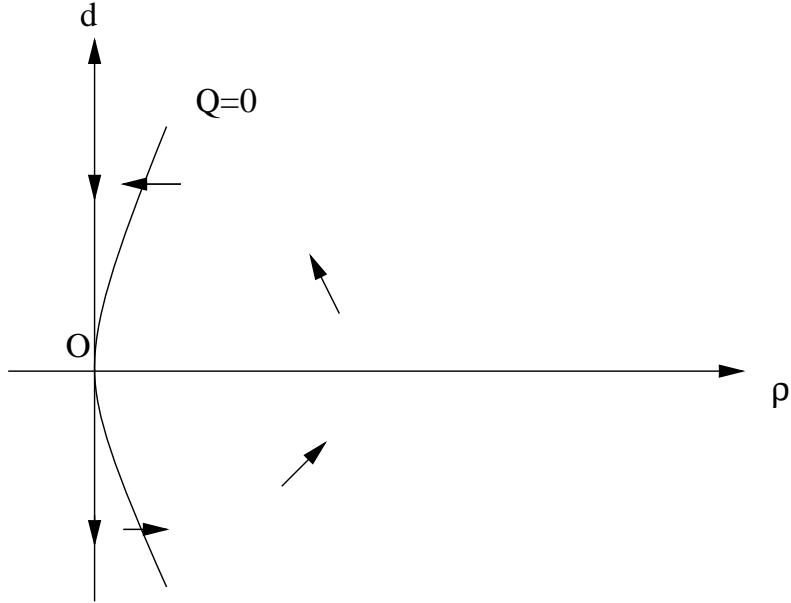

FIGURE 2.1. Zero level set of $d' = 0$. Complex spectral gap

The trajectory on the plane $\rho < 0$ does not affect the solution behavior in the region $\rho > 0$ since $\rho = 0$ is an invariant set, governed by

$$\rho \equiv 0, \quad d' = -\frac{d^2}{2} \rightarrow d(t) = \frac{d_0}{1 + \frac{d_0}{2}t}.$$

Note that $(0,0)$ is the only critical point of the autonomous ODE system (2.6), (2.7) on the right half phase plane, and that the vector field in $\{(\rho,d), \quad Q < 0, d \geq 0\}$ is converging to the critical point $(0,0)$. It follows that for global smoothness it suffices to control the divergence $d$ from below in the region $\{(\rho,d), \quad Q < 0, d < 0, \rho > 0\}$, and to control the density from above in the region $\{(\rho,d), Q > 0\}$.

For the former case we have, recalling that $\beta < 0$,

$$\left(\frac{d}{\rho}\right)' = k + \frac{d^2}{2\rho} - \frac{\beta\rho}{2} \geq k,$$

and its integration along particle path gives

$$d \geq \left(kt + \frac{d_0}{\rho_0}\right)\rho.$$

This shows that the divergence $d$ is bounded from below, and in particular, it becomes positive for large time since the density is positive. To the upper-bound for $\rho$ in the region $Q > 0$, where $d(t) \geq d_0(\alpha)$, we substitute this estimate into the mass equation yielding

$$\rho' \leq -\rho d_0(\alpha).$$

This clearly gives the upper-bound for the density $\rho \leq \rho_0(\alpha)e^{-d_0 t}$.

2.2. **Real Spectral Gap.** When $\beta \geq 0$ the initial spectral gap is real, and there are two cases to be considered:

**Subcase 1:** $\beta = 0$ when the eigenvalues are equal, i.e. $\lambda_1(0) = \lambda_2(0)$.

In this case the zero level set $d' = 0$ becomes a parabola passing through the only critical point $(0,0)$, see Figure 2.2, and one can repeat similar augments to our phase plane analysis in the previous case of distinct real roots. Note that the global invariant (2.8) becomes

$$\frac{d^2}{\rho} + 2kln\rho = Const.$$

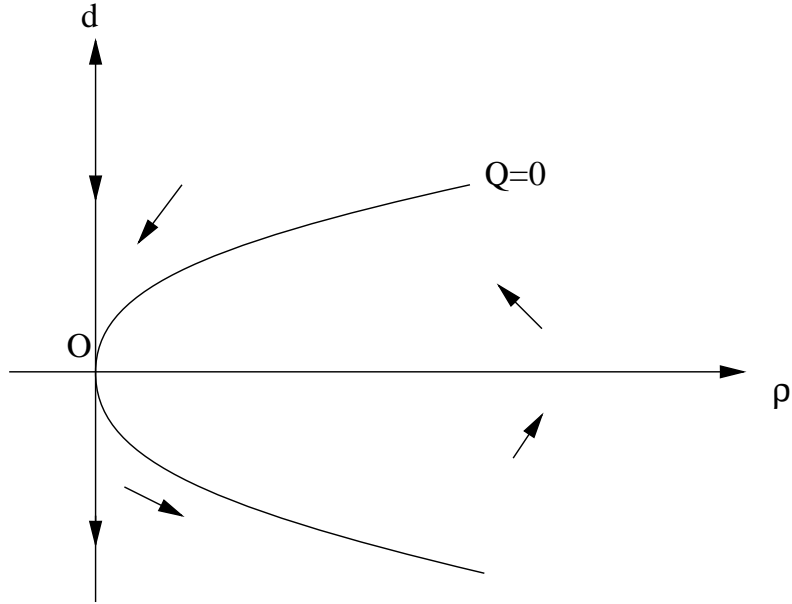

FIGURE 2.2. Zero level set $d' = 0$. Real spectral gap.

**Subcase 2:** $\beta > 0$ when the eigenvalues are initially real.

**Lemma 2.3.** *If eigenvalues of $\nabla U_0$ are real, then the solution of 2-D REP remains smooth for all time if and only if*

$$\lambda_1(0) + \lambda_2(0) \geq g(\rho_0),$$

*where*

$$g(\rho) := \operatorname{sgn}(\rho - \frac{2k}{\beta})\sqrt{\rho \times \left(F(\rho) - F\left(\frac{2k}{\beta}\right)\right)}, \quad F(\rho) = \beta\rho - 2kln\rho.$$

*Proof.* The system (2.6)-(2.7) has two critical points on the phase plane: $O(0,0)$ and $A(\frac{2k}{\beta}, 0)$, see Figure 2.3.

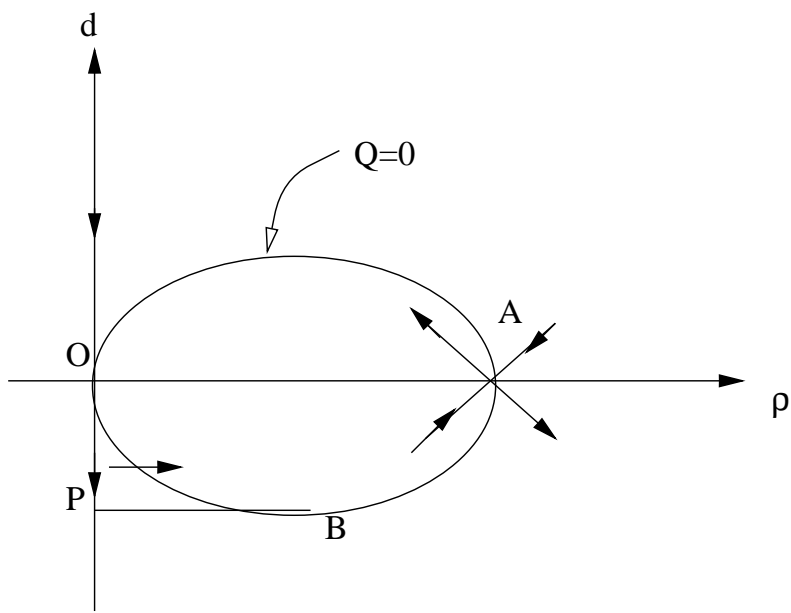

FIGURE 2.3. Critical points in $\rho - d$ plane. Real spectral gap

The coefficient matrix of the linearized system around $(\rho^*, d^*)$ is
$$L(\rho^*, d^*) = \begin{pmatrix} -d^* & -\rho^* \\ k - \beta\rho^* & -d^* \end{pmatrix}.$$
A simple calculation gives the eigenvalues of $L$,
$$\lambda_\pm = -d^* \pm \sqrt{\rho^*(\beta\rho^* - k)}.$$
At $(0,0)$, we have $\lambda_1 = \lambda_2 = 0$ and hence $(0,0)$ is a non hyperbolic critical point. Another critical point, $A(\frac{2k}{\beta}, 0)$, is a saddle since $\lambda_{1,2} = \pm\sqrt{\frac{2}{\beta}}k$. We shall use the above facts to construct the critical threshold via the phase plane analysis.

Assume the seperatrix enters (leaves) $A$ along the line $d = s(\rho - \frac{2k}{\beta})$. Upon substitution into the linearized system around $A$, i.e.,
$$\rho' = -\frac{2k}{\beta}d, \quad d' = -k\left(\rho - \frac{2k}{\beta}\right),$$
one can obtain
$$s = \pm\sqrt{\frac{\beta}{2}}.$$
Thus, two seperatrixes leave/enter $A$ along the directions
$$\theta_1 = -arctg\sqrt{\frac{\beta}{2}} \quad \text{and} \quad \theta_2 = arctg\sqrt{\frac{\beta}{2}}.$$
In the phase plane the zero level set of $d' = 0$ is an ellipse, see Figure 2.3,
$$d^2 + \beta\left(\rho - \frac{k}{\beta}\right)^2 = \frac{k^2}{\beta}.$$

Let $\gamma_s(A)$ be the portion of the stable manifold of the system coming into $A$ from $\{d < 0\}$. In order to prove the existence of a critical threshold it suffices to show that $\gamma_s(A)$ can only come from $O$. Let $B$ be the lowest point of the ellipse with coordinates $\left(\frac{k}{\beta}, -\frac{k}{\sqrt{\beta}}\right)$ and $\overline{PB}$ be a horizontal line intersecting with $\rho = 0$ at $\left(0, -\frac{k}{\sqrt{\beta}}\right)$. According to the vector field inside the ellipse we see that the trajectory $\gamma_s(A)$ can only come from the area $OPB$ by crossing the curve $\overline{OB}$. Note that the vector field on $\overline{PB}$ is going outside $OPB$ and $\rho = 0$ is invariant. Thus all trajectories in the area $OPB$ originate from $O$. Therefore $\gamma_s(A)$ can only originate from $O$ (as $t \to -\infty$) and becomes a portion of one unstable manifold of $O$. By symmetry we can show that the unstable manifold of the system issued from $A$ entering $\{d > 0\}$ will end through the portion $\{d > 0\}$ at $O$.

Thus the critical curve $g : \mathbb{R}^+ \to \mathbb{R}$ is the one defined as

$$\{(\rho, d), \quad d = g(\rho)\} = \gamma_s(A).$$

In order to have a precise formula for $g$ we need to use the global invariant (2.8), i.e.,

$$\frac{d^2 - \beta\rho}{\rho} + 2kln\rho = Const.$$

Thus all trajectories can be expressed as

$$d^2 = \rho[C(\alpha) + F(\rho)]$$

with

$$C(\alpha) = \frac{d_0^2}{\rho_0} - F(\rho_0), \quad F(\rho) := \beta\rho - 2k \ln \rho.$$

Note that $F(\rho)$ is a convex function and $\min_{\rho>0} F = F(\frac{2k}{\beta}) = 2k[1 - ln(2k/\beta)]$. Due to the symmetry the homoclinic connection is possible when the trajectory pass $(\rho_0, 0)$ with $\rho_0 \leq \frac{2k}{\beta}$ and converging to $(0,0)$ as $t \to \pm\infty$, i.e. the initial data must satisfy $0 < \rho_0 < \frac{2k}{\beta}$ and

$$C(\alpha) \leq -F\left(\frac{2k}{\beta}\right), \quad i.e. \quad \frac{d_0^2}{\rho_0} \leq F(\rho_0) - F\left(\frac{2k}{\beta}\right).$$

The seperatrices passing through $(\frac{2k}{\beta}, 0)$ correspond to $C(\alpha) = -F(\frac{2k}{\beta})$. The stable manifold $\gamma_s(A)$ can be written as $d = g(\rho)$ for $0 \leq \rho < \infty$, where

$$g(\rho) = \text{sgn}\left(\rho - \frac{2k}{\beta}\right)\sqrt{\rho\left(F(\rho) - F\left(\frac{2k}{\beta}\right)\right)}.$$

It remains to prove that the initial data satisfying $d_0 < g(\rho_0)$ always lead to finite time breakdown.

First, in the region $\left\{(\rho, d), \quad d < -\sqrt{\rho(F(\rho) - F(\frac{2k}{\beta}))}\right\}$, there must exist a finite time $T_1 > 0$ such that $\rho(T_1) > \frac{2k}{\beta}$ for $\rho_0 \leq \frac{2k}{\beta}$ ( take $T_1 = 0$ for $\rho_0 > \frac{2k}{\beta}$) since $\rho' > 0$. Therefore $\rho(t) \geq \rho(T_1)$ for $t \geq T_1$ and

$$\rho' = -\rho d \geq \rho\sqrt{\rho\left(\rho(T_1) - F\left(\frac{2k}{\beta}\right)\right)} \quad \text{for} \quad t \geq T_1.$$

Integration over $[T_1, t]$ gives

$$\sqrt{\rho(t)} \geq \frac{\sqrt{2\rho(T_1)}}{2 - (t - T_1)\sqrt{\rho(T_1)(F(\rho(T_1)) - F(\frac{2k}{\beta}))}}, \quad t \geq T_1.$$

Thus the solution must becomes unbounded before the time

$$T_1 + \frac{2}{\sqrt{\rho(T_1)(F(\rho(T_1)) - F(\frac{2k}{\beta}))}}.$$

Second, we consider the trajectories in the region

$$\left\{ (\rho, d), \quad \rho > \frac{2k}{\beta}, \quad |d| < \sqrt{\rho\left(F(\rho) - F\left(\frac{2k}{\beta}\right)\right)} \right\}.$$

Note that at finite time the trajectory must enter the subregion $\{(\rho, d), \quad d < 0\}$ through the left point $(\rho^*, 0)$ identified as

$$d^2 = \rho[F(\rho) - F(\rho^*)], \quad \rho \geq \rho^* > \frac{2k}{\beta}.$$

This combined with the Riccati-type inequality

$$d' < -d^2/2$$

ensures the breakdown at finite time. This completes the confirmation of the curve $d = g(\rho)$ as a critical threshold. □

**Proof of Theorem 1.1:** It suffices to summarize the above cases with

$$\beta = \frac{\Gamma_0}{\rho_0^2}$$

being taken into account. Clearly the case $\beta < 0$ and $\beta = 0$ correspond to the set

$$\{(\rho_0, M_0), \quad \Gamma_0 \leq 0\}.$$

For $\beta > 0$, i.e., $\Gamma_0 > 0$ we rewrite the critical threshold as

$$d_0 = \text{sgn}(\rho_0 - \frac{2k}{\beta})\sqrt{\rho_0 \left(F(\rho_0) - F\left(\frac{2k}{\beta}\right)\right)}$$

$$= \text{sgn}(\Gamma_0 - 2k\rho_0)\sqrt{\Gamma_0 - 2k\rho_0 + 2k\rho_0 ln\left(\frac{2k\rho_0}{\Gamma_0}\right)},$$

where we have used the relation $F(\rho) = \beta\rho - 2kln\rho$ and

$$F(\rho_0) = \frac{\Gamma_0}{\rho_0} - 2kln\rho_0,$$

$$F\left(\frac{2k}{\beta}\right) = 2k - 2kln\left(\frac{2k\rho_0^2}{\Gamma_0}\right).$$

This completes the proof of Theorem 1.1. □

## 3. 2D REP with non-zero background

This section is devoted to the study of the REP with nonzero background $c > 0$, for which the velocity gradient tensor, $M = \nabla U$, solves

$$\text{(3.1)} \quad \frac{d}{dt}M + M^2 = \frac{k}{2}[\rho - c],$$

$$\text{(3.2)} \quad \frac{d}{dt}\rho + \rho\, tr M = 0.$$

Again, using the Spectral Dynamics Lemma presented in [17] the spectral dynamics of $M$ is governed by

$$\lambda_1' + \lambda_1^2 = \frac{k}{2}(\rho - c), \quad ' := \frac{d}{dt}$$

$$\lambda_2' + \lambda_2^2 = \frac{k}{2}(\rho - c),$$

$$\rho' + \rho(\lambda_1 + \lambda_2) = 0.$$

As in the zero background case, the difference $\eta := \lambda_2 - \lambda_1$ is proportional to the density along the particle path in the sense that

$$\frac{\eta(t)}{\rho(t)} = \frac{\eta_0(\alpha)}{\rho_0(\alpha)}, \quad \alpha \in \mathbb{R}^2.$$

Further manipulation gives a closed system

$$\text{(3.3)} \quad \rho' = -\rho d$$

$$\text{(3.4)} \quad d' = k(\rho - c) - \frac{d^2 + \beta\rho^2}{2} =: Q, \quad \beta := \eta_0^2/\rho_0^2.$$

Once again the dynamics of (3.3), (3.4) is influenced by the choice of $\beta$. We proceed to discuss the solution behavior of (3.3), (3.4) by distinguishing two cases:
(1) $\beta < 0$, the spectral gap is complex;
(2) $\beta \geq 0$, the spectral gap is real.

### 3.1. Complex Spectral Gap.
We first discuss the case $\beta < 0$, which corresponds to the case that the eigenvalues are initially complex.

**Lemma 3.1.** *Assume that the eigenvalues are initially complex with $Im(\lambda_i(0)) \neq 0$. Then the solution of (3.3), (3.4) remains smooth for all time. Moreover, there is a global invariant in time, given by*

$$V(\rho, d) = \rho^{-1}\left[d^2 - \beta\rho^2 + 2k\rho ln\left(\frac{\rho}{2c}\right) + 2ck\right],$$

*Proof.* A straightforward computation yields $\dot{V} = 0$ along the 2D REP solutions, which implies that the curves $V = Const$ are invariants of the flow. As before, for negative $\beta$'s we have

$$d^2 \leq \max_{\rho > 0}\left\{Const.\rho - 2k\rho ln\left(\frac{\rho}{2c}\right) - 2ck + \beta\rho^2\right\} \leq C_1^2$$

and the bounds of $d$ (and hence of $\rho$) follow. □

*Remark.*

1. As before, a more detailed information is available in this case by a phase plane analysis. If the eigenvalues are initially complex, then one has $\beta = \eta_0^2/\rho_0^2 < 0$. The zero level set of $d' = Q = 0$ becomes a hyperbola, see Figure 3.1.

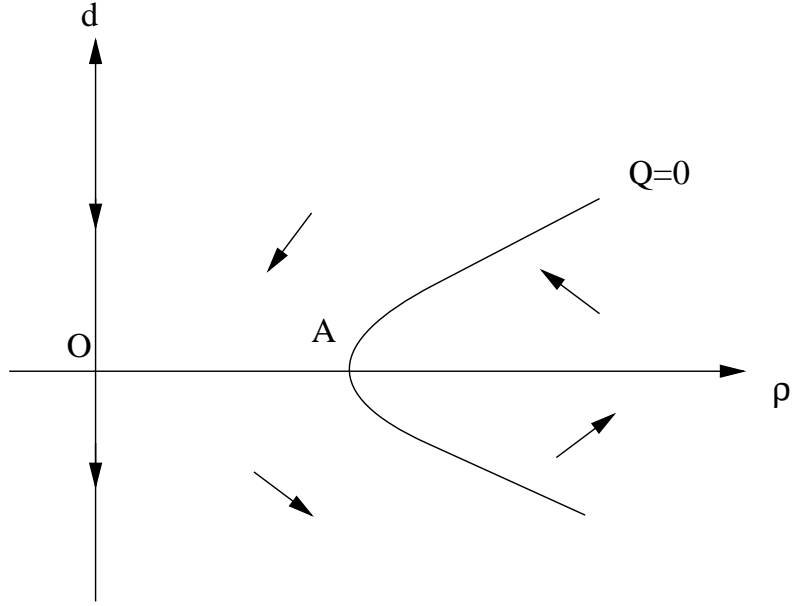

FIGURE 3.1. Zero level set $d' = 0$. Complex spectral gap with nonzero background.

The intersection of its right branch with $d = 0$ is the rest point $A = (\rho^*, 0)$ of the system, where
$$\rho^* = \frac{k}{\beta} + \sqrt{\frac{k^2}{\beta^2} - \frac{2ck}{\beta}}.$$
The coefficient matrix of the linearization around $(\rho^*, 0)$ is
$$L(\rho^*, 0) = \begin{pmatrix} 0 & -\rho^* \\ k - \beta\rho^* & 0 \end{pmatrix}.$$
Its eigenvalues satisfy
$$\lambda^2 = \rho^*(\beta\rho^* - k) = -\rho^*\sqrt{k^2 - 2ck\beta} < 0.$$
Hence such critical point is a non hyperbolic equilibrium. The nonlinear effect plays essential roles in the solution behavior. In order to locate the possible critical threshold, we first study the solution around $(\rho^*, 0)$. Set $n = \rho - \rho^*$, we then have

(3.5)  $$n' = -\rho^* d - nd,$$

(3.6)  $$d' = \sqrt{k^2 - 2ck\beta}\, n - \frac{d^2}{2} - \frac{\beta}{2} n^2.$$

It is easy to see that the flow governed by the linear part stays on the ellipse
$$\sqrt{k^2 - 2ck\beta}\, n^2 + \rho^* d^2 = Const.$$

In order to capture the dynamics of the nonlinear system in the neighborhood of the critical point $(n, d) = (0, 0)$, we employ the polar coordinates of the form

$$n = \frac{r\cos\theta}{(k^2 - 2ck\beta)^{1/4}},$$

$$d = \frac{-r\sin\theta}{\sqrt{\rho^*}}.$$

Careful calculation with these polar coordinates, yields that (3.5(-(3.6) recast into the form

(3.7) $$r' = R(r, \theta),$$
(3.8) $$\theta' = -\sqrt{\rho^*}(k^2 - 2ck\beta)^{1/4} + \Theta(r, \theta),$$

where

$$R(r, \theta) = \frac{r^2 \sin\theta}{2\sqrt{\rho^*}} \left[1 + \frac{k\cos^2\theta}{\sqrt{k^2 - 2ck\beta}}\right]$$

$$\Theta(r, \theta) = -\frac{r\cos\theta}{2\sqrt{\rho^*}\sqrt{k^2 - 2ck\beta}} \left[\sqrt{k^2 - 2ck\beta}\sin^2\theta - \beta\rho^*\cos^2\theta\right].$$

When $r$ is sufficiently small $\theta'$ is strictly negative. The pleasant implication of this is that the orbits of system (3.5), (3.6) spiral monotonically in $\theta$ around the $(\rho^*, 0)$. But the even power of $r^2$ does not tells the stability property of the critical point.

Observe that if $(n(t), d(t))$ is a solution, so is $(n(-t), -d(-t))$. Such symmetry implies that there is a center in the neighborhood of $(\rho^*, 0)$.

In order to clarify the global behavior of the flow around such center, we appeal to the global invariant

$$V(\rho, d) = \rho^{-1} \left[d^2 - \beta\rho^2 + 2k\rho \ln\left(\frac{\rho}{2c}\right) + 2ck\right].$$

we claim that $V$ is positive definite which serves as a (majoration of) Lyapunov functional. to this end, we consider the function $H(\rho) := -\beta\rho^2 + 2k\rho \ln\left(\frac{\rho}{2c}\right) + 2ck$, which is convex and takes its minimum at $\rho_{\min}$, satisfying

$$\ln\left(\frac{\rho_{\min}}{2c}\right) = -1 + \frac{\beta}{k}\rho_{\min}.$$

Observe that since $\beta < 0$ the function

$$h(\rho) := 1 - \frac{\beta}{k}\rho + \ln\left(\frac{\rho}{2c}\right)$$

is an increasing function in $\rho > 0$ and $h(\rho_{\min}) = 0$, which when combined with the fact that $h(\rho^*) \geq 0$ verifies that

$$0 < \rho_{\min} \leq \rho^*, \qquad \rho^* := \beta^{-1}[k - \sqrt{k^2 - 2ck\beta}].$$

Indeed, for $\rho^*$ we have

$$h(\rho^*) = 1 - \frac{\beta}{k}\rho^* + \ln\left(\frac{\rho^*}{2c}\right) = \sqrt{1 - 2ck^{-1}\beta} - \ln\left(1 + \sqrt{1 - 2ck^{-1}\beta}\right) \geq 0.$$

Therefore $H(\rho)$ is nonnegative since

$$H(\rho_{\min}) = \beta\rho_{\min}^2 - 2k\rho_{\min} + 2ck = \beta(\rho_{\min} - \rho^*)(\rho_{\min} - \bar{\rho}^*) \geq 0,$$

where $\bar{\rho}^* = \beta^{-1}[k + \sqrt{k^2 - 2ck\beta}]$.

The invariant curves, $V = Const.$, represent of course, the bounded periodic orbits containing $(\rho^*, 0)$.

### 3.2. Real Spectral Gap.

We divide the region $\beta \in [0, \infty)$ into subregions depending on the number of critical points on the phase plane, and then study the solution behavior with $\beta$ in each sub-region. The solution behavior depends strongly on the number of critical points and their stability property.

Let $(\rho^*, d^*)$ be a critical point of the system, then the coefficient matrix of the linearization around $(\rho^*, d^*)$ reads

$$L(\rho^*, d^*) = \begin{pmatrix} -d^* & -\rho^* \\ k - \beta\rho^* & -d^* \end{pmatrix}.$$

Its eigenvalues are given by

$$\lambda = -d^* \pm \sqrt{\rho^*(\beta\rho^* - k)}. \tag{3.9}$$

We now discuss subcases distinguished by the number and type of the critical points as $\beta$ changes.

- $\beta = 0$, then the zero level set $d' = Q = 0$ is a parabola, $d^2 = 2k(\rho - c)$, intersecting with $d = 0$ at $(\rho^*, d^*) = (c, 0)$. From (3.9) we see that at this point the eigenvalues of $L$ are $\lambda = \pm\sqrt{ck}i$, pure imaginary number, the critical point $(c, 0)$ is non-hyperbolic. The stability property of this critical point has to be determined by taking into account of the nonlinear effect.
- $0 < \beta < \frac{k}{2c}$. The zero level set $Q = 0$ is an ellipse, located on the right half plane $\rho > 0$. There are two critical points $(\rho^*, d^*) = (\rho^*, 0)$ with

$$\rho^* = \frac{k}{\beta} \pm \sqrt{\frac{k^2}{\beta^2} - \frac{2kc}{\beta}}.$$

The associated eigenvalues of $L$ are

$$\lambda(\rho_1^*) = \pm\sqrt{\rho_1^*\sqrt{k^2 - 2ck\beta}}i, \quad \lambda(\rho_2^*) = \pm\sqrt{\rho_2^*\sqrt{k^2 - 2ck\beta}}.$$

Therefore $(\rho_1^*, 0)$ is a center of the linearized system and $(\rho_2^*, 0)$ is a saddle, see Figure 3.2.

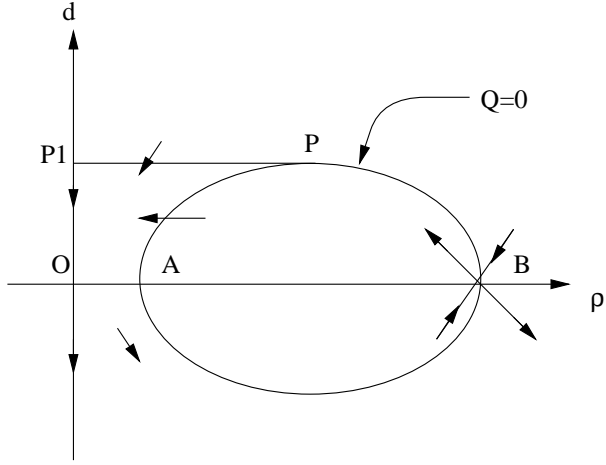

FIGURE 3.2. Critical points in $\rho - d$ plane. Real spectral gap with nonzero background.

Possible bifurcation as $\beta$ changes from 0 to $\frac{k}{2c}$ may be responsible for the complicated solution structure in this regime.
- $\beta = \frac{k}{2c}$. The zero level set $Q = 0$, i.e.,
$$d^2 + \beta(\rho - \frac{k}{\beta})^2 = 0,$$
degenerates to a single point $(\rho^*, d^*) = (\frac{k}{\beta}, 0)$, the only critical point with zero eigenvalues.
- $\beta > \frac{k}{2c}$. In this case
$$Q = -\frac{1}{2}\left[d^2 + \beta\left(\rho - \frac{k}{\beta}\right)^2 + 2kc - \frac{k^2}{\beta}\right] \leq \frac{k^2}{2\beta} - kc < 0.$$

There is no critical point at all in the finite phase plane.

The solution behavior distinguished by above cases is given in the following lemmata.

**Lemma 3.2.** *If $\lambda_1(0) = \lambda_2(0)$. Then the solution of (3.3), (3.4) remains smooth for all time, indicated by the global invariant*

$$\frac{d^2 + 2ck}{\rho} + 2k\ln\rho = Const. \tag{3.10}$$

*Proof.* The assumption amounts to $\beta = 0$. As discussed above, $(c, 0)$ is the only critical point and the center of the linearized system. In order to find the global invariant we set $R(t) := k(\rho - c)^2 + cd^2$. Along the trajectory $\frac{dx}{dt} = U(x, t)$

$$\frac{d}{dt}R(t) = 2k(\rho - c)\rho' + 2cdd' = -d[k(\rho - c)^2 + R(t)]. \tag{3.11}$$

¿From the mass equation it follows that

$$d = -\frac{\rho'}{\rho},$$

which when inserted into the relation (3.11) gives

$$\frac{dR}{d\rho} = \frac{k(\rho - c)^2}{\rho} + \frac{R}{\rho}.$$

Integration gives

$$\frac{R}{\rho} + \frac{kc^2}{\rho} + 2ckln\rho - k\rho = Const,$$

which leads us to the global invariant as asserted in (3.10). This global invariant is compact and ensures that both divergence $d$ and the density $\rho$ remains bounded as time evolves. □

We leave the case $0 < \beta < \frac{k}{2c}$ later and study the critical case $\beta = \frac{k}{2c}$.

**Lemma 3.3.** *If $\lambda_2(0) - \lambda_1(0) = \sqrt{\frac{k}{2c}}\rho_0$. Then the solution of (3.3), (3.4) always develops finite time breakdown unless the initial data lies in the set*

$$\left\{ (\rho, d) \in \mathsf{IR}^+ \times \mathsf{IR}^+, \quad \frac{d^2 + 2ck}{\rho} - \frac{k}{2c}\rho + 2kln\rho = 2kln(2c) \right\}.$$

*Proof.* The given assumption is equivalent to the case $\beta = \frac{k}{2c}$. In this case the divergence always decrease except for at the critical point $(2c, 0)$ since

$$d' = -\frac{d^2}{2} - \frac{k(\rho - 2c)^2}{4c} \leq 0.$$

In order to clarify the solution behavior we proceed to obtain the global invariant. Set $q := d^2$, one then has

$$\frac{dq}{d\rho} = \frac{2dd'}{\rho'} = \frac{q + \beta\rho^2 - 2k(\rho - c)}{\rho}.$$

Solving this equation we obtain

$$\frac{q}{\rho} = \beta\rho - 2kln\rho - \frac{2ck}{\rho} + Const.$$

Therefore we come up with a global invariant

(3.12) $$\frac{d^2 + 2ck}{\rho} - \beta\rho + 2kln\rho = Const.$$

The only trajectory converging to the critical point is realized by a half trajectory converging to $(2c, 0)$ from the first quadrant. For all other trajectories not passing the critical point $(2c, 0)$, the rate $d'$ is strictly negative. The divergence will become negative at finite time even if it is initially positive, which when combined with the Riccati-type inequality $d' \leq -d^2/2$ confirms the finite time breakdown. □

We now look at the case $\beta > \frac{k}{2c}$.

**Lemma 3.4.** *Assume that the eigenvalues are initially real and $|\lambda_2(0) - \lambda_1(0)| > \sqrt{\frac{k}{2c}}\rho_0(\alpha)$. Then the solution of (3.3), (3.4) always develops finite time breakdown.*

*Proof.* The given assumption is nothing but the inequality $\beta > \frac{k}{2c}$. Note that there is no critical point in the finite phase plane, actually $Q$ remains negative for all time. The solution must develops breakdown in finite time. In fact from

$$(3.13) \qquad d' = -\frac{d^2}{2} - \frac{\beta}{2}\left(\rho - \frac{k}{\beta}\right)^2 - \frac{kc}{\beta}\left(\beta - \frac{k}{2c}\right),$$

we find that

$$d' \leq -\delta \quad \text{with} \quad \delta := \frac{kc}{\beta}\left(\beta - \frac{k}{2c}\right) > 0.$$

This ensures that $d$ must become negative beyond a finite time $T_0$, say $T_0 > \max\{\frac{d_0}{\delta}, 0\}$. The $d-$ equation (3.13) also gives

$$d' \leq -d^2/2,$$

whose integration over $[T_0, t]$ leads to

$$d(t) \leq \frac{d(T_0)}{1 - \frac{1}{2}d(T_0)(t - T_0)}.$$

Hence the solution must break down at a finite time before $T_0 - \frac{2}{d(T_0)}$. □

Finally we conclude this subsection by discussing the delicate case $0 < \beta < \frac{k}{2c}$. Set

$$G(\rho, \rho^*, \beta) := \beta(\rho - \rho^*) - 2k\ln\left(\frac{\rho}{\rho^*}\right) - \frac{2ck}{\rho} + \frac{2ck}{\rho^*}$$

with $\rho^* = \beta^{-1}[k \pm \sqrt{k^2 - 2ck\beta}]$ being the $\rho$-coordinate of the intersection point of the trajectory with $\rho$ axis.

**Lemma 3.5.** *Assume that the real eigenvalues satisfy $0 < |\lambda_2(0) - \lambda_1(0)| < \sqrt{\frac{k}{2c}}\rho_0(\alpha)$. Then for any $\beta \in (0, \frac{k}{2c})$, the solutions of (3.3), (3.4) remains smooth for all time if and only if*

$$|\lambda_1(0) + \lambda_2(0)| \leq \sqrt{\rho_0 G(\rho_0, \rho_2^*, \beta_0)} \quad \text{for} \quad \rho_0 \leq \rho_2^*$$

*and*

$$\lambda_1(0) + \lambda_2(0) = \sqrt{\rho_0 G(\rho_0, \rho_2*, \beta_0)} \quad \text{for} \quad \rho_0 > \rho_2^*.$$

*Proof.* The assumption tells us that $\beta < \frac{k}{2c}$. In this case there are two critical points in the phase plane $A = (\rho_1^*, 0)$ and $B = (\rho_2^*, 0)$, see Figure 3.2. $B$ is a saddle, whose two manifolds pass enclosing the critical point $A$, which is a center of the linearized system. Let $W_s(B)$ denote the stable manifold coming from the region $\{(\rho, d), \rho < \rho_2^*, d < 0\}$ and $W_u(B)$ the unstable manifold entering into $\{(\rho, d), \rho < \rho_2^*, d > 0\}$. To prove the results stated in the theorem it suffices to show that for any $\beta \in (0, \frac{k}{2c})$ such that

$$W_u(B) \cap W_s(B)$$

is not empty, i.e. there always exists a saddle connection (homoclinic orbit).

This follows from the continuity argument supported by the following facts:

(1) Both $W_u(B)$ and $W_s(B)$ pass through the segment $OA$ with flow going downward since $d' < 0$ and $\rho' = 0$ on $OA$, see Figure 3.2

The level curve $d' = 0$ is an ellipse with upper vortex $P$ located at $(\frac{k}{\beta}, 0)$. Let $P_1$ denote the intersection of tangent line of the ellipse $Q = 0$ through $P$ with the axis $\rho = 0$. The vector field inside the ellipse shows that $W_u(B)$ must escape the ellipse from the curve $\overline{PA}$. Note that the trajectories on $\overline{PP_1}$ and $\overline{AP}$ enters into the region $PAOP_1$, and the axis $\rho = 0$ is an invariant set. These facts ensure that $W_u(B)$ must enter the region $d < 0$ through $OA$. Similarly we can show that $W_s(B)$ for $\rho \leq \rho_2^*$ must enter the region $d > 0$ through $OA$.

(2) As $\beta$ increase in $(0, \frac{k}{2c})$, the point $W_u(B) \cap OA$ moves to the right and the point $W_s(B) \cap OA$ moves to the left.

We prove the claim for $W_u(B) \cap OA$, and the case for $W_s(B) \cap OA$ follows similarly. The claim follows from the following two observations:

(i) The slope of the unstable manifold $W_u(B, \beta)$ at $(\rho_2^*, 0)$ is $\partial_\rho d|_{\rho=\rho_2^*} = \lambda_-(\rho_2^*, \beta)$ and the eigenvalue $\lambda_-(\rho_2^*, \beta)$ is increasing in $\beta$; Indeed

$$\frac{d\lambda_-(\rho_2^*, \beta)}{d\beta} = -\frac{\beta^2}{k\lambda_-(\rho_2^*, \beta)}\left\{k + \sqrt{k^2 - 2ck\beta} + \frac{ck\beta}{\sqrt{k^2 - 2ck\beta}}\right\} > 0.$$

(ii) $W_u(B, \beta_1)$ does not intersect with $W_u(B, \beta_2)$ for $\beta_1 \neq \beta_2$;

As done previously we can find the global invariant of the system

$$d^2 = \rho\left[\beta\rho - 2kln\rho - \frac{2ck}{\rho} + Const\right],$$

from which the left branch of the unstable manifold of $B$ can be explicitly expressed as

$$d = \sqrt{\rho\left[\beta(\rho - \rho_2^*) - 2kln(\rho/\rho_2^*) - \frac{2ck}{\rho} + \frac{2ck}{\rho_2^*}\right]}, \quad 0 < \rho < \rho_2^*.$$

A careful calculation gives

$$\frac{\partial d}{\partial \beta} = \frac{\rho}{2d}(\rho - \rho_2^*) < 0,$$

which ensures the claim (ii).

(3) Let $\rho_u(\beta)$ be the $\rho$ coordinate of the point $W_u(B, \beta) \cap OA$ and $\rho_s(\beta)$ be $\rho$-coordinate of the point $W_s(B, \beta) \cap OA$. We claim

$$\lim_{\beta \to 0+} \rho_u(\beta) < \lim_{\beta \to \frac{k}{2c}-} \rho_s(\beta).$$

In fact from the expression of seperatrices

$$d^2 = \rho G(\rho, \rho_2^*, \beta), \quad \rho < \rho^*,$$

we find that the $\rho$-coordinate of points $W_{u/s}(B) \cap OA$ satisfy

$$G(\rho, \rho_2^*, \beta) \equiv 0.$$

Note that

$$\frac{\partial G}{\partial \rho} = \frac{\beta}{\rho^2}(\rho - \rho_1^*)(\rho - \rho_2^*), \quad \frac{\partial G}{\partial \beta} = \rho - \rho_2^*.$$

Thus we have for $0 < \rho < \rho_1^*$

$$\frac{\partial \rho}{\partial \beta} = -\frac{\partial G}{\partial \beta}\bigg/\frac{\partial G}{\partial \rho} = \frac{\rho^2}{\beta}(\rho_1^* - \rho) > 0.$$

This confirms the above assertion.

Combining the above observations we conclude that there exists a $\beta_0 \in (0, \frac{k}{2c})$ for which a saddle connection exists. It remains to show that as $\beta$ changes in the region $(0, \frac{k}{2c})$ the above saddle connection is preserved. Observe that if $(\rho(t), d(t))$ is a solution, so is $(\rho(-t), -d(-t))$. Such symmetry prevents the occurrence of the possible bifurcation when $\beta$ changes.

Using the nonlinear terms in the equation and the vector field we can show for the initial data outside the closed curve — saddle connection — the solution always develops finite time breakdown, details are omitted. □

**Proof of Theorem 1.2**: Summarizing the results stated in the above lemmas we see that the case $\beta < 0$ and $\beta = 0$ corresponds to the set

$$S_1 = \left\{ (\rho_0, d_0, \Gamma_0) \;\bigg|\; \Gamma_0 \leq 0 \quad \text{and} \quad \begin{cases} d_0 \geq 0 & \text{if } \rho_0 = 0 \\ d_0 \text{ arbitrary} & \text{if } \rho_0 > 0 \end{cases} \right\},$$

since $\Gamma_0 = \beta \rho_0^2$. The case $0 < \beta < \frac{k}{2c}$ corresponds to $0 < \Gamma_0 < \frac{k}{2c} \rho_0^2$, the divergence is required to satisfy the critical threshold condition

$$|d_0| \leq \sqrt{\rho_0 G(\rho_0, \rho_2^*, \beta)}, \quad 0 < \rho_0 < \rho_2^*$$

and $d_0 = \sqrt{\rho_0 G(\rho_0, \rho_2^*, \beta)}$ for $\rho_0 \geq \rho_2^*$. Using $\Gamma_0 = \beta \rho_0^2$ and

$$\rho_2^* = \beta^{-1}[k + \sqrt{k^2 - 2ck\beta}] = \frac{2ck}{k^2 - \sqrt{k^2 - 2ck\beta}} = \frac{2c\rho_0}{\rho_0 - \sqrt{\rho_0^2 - \frac{2c}{k}\Gamma_0}},$$

one has

$$\rho_0 G(\rho_0, \rho_2^*, \beta) = \rho_0 \left[ \beta(\rho_0 - \rho_2^*) - 2k \ln\left(\frac{\rho_0}{\rho_2^*}\right) - \frac{2ck}{\rho_0} + \frac{2ck}{\rho_2^*} \right]$$

$$= \Gamma_0 \left( 1 - \frac{2c}{\rho_0 - \sqrt{\rho_0^2 - \frac{2c\Gamma_0}{k}}} \right) - 2k\rho_0 \ln\left( \frac{\rho_0 - \sqrt{\rho_0^2 - \frac{2c\Gamma_0}{k}}}{2c} \right)$$

$$- 2ck + k\left( \rho_0 - \sqrt{\rho_0^2 - \frac{2c\Gamma_0}{k}} \right)$$

$$= \Gamma_0 - 2ck - 2k\sqrt{\rho_0^2 - \frac{2c\Gamma_0}{k}} - 2k\rho_0 \ln\left( \frac{\rho_0 - \sqrt{\rho_0^2 - \frac{2c\Gamma_0}{k}}}{2c} \right),$$

which leads to the critical threshold described by the set $S_2$. The set $S_3$ can be figured out in a similar manner.

□

## ACKNOWLEDGMENTS

Research was supported in part by ONR Grant No. N00014-91-J-1076 (ET) and by NSF grant #DMS01-07917 (ET, HL).

UCLA, Mathematics Department, Los Angeles, CA 90095-1555.
*E-mail address*: `hliu@math.ucla.edu`

UCLA, Mathematics Department, Los Angeles, CA 90095-1555.
*E-mail address*: `tadmor@math.ucla.edu`